\documentclass[a4paper,10pt]{amsart}
\usepackage{amsmath}
\usepackage{amssymb}
\usepackage{amsfonts}
\usepackage{amscd}
\usepackage{latexsym,graphicx,caption}
\input xy
\xyoption{all}
\usepackage{verbatim}
\newcommand{\ti}[1]{\mbox{$\tilde{#1} $}}

\newcommand{\lo}{\longrightarrow}
\newcommand{\llo}{\longleftarrow}

\newcommand{\diam}{{\hfill \nobreak} $\Box$}
\newcommand{\noi}{\noindent}
\newcommand{\lspace}{\vspace{0.45cm}}
\newcommand{\sspace}{\vspace{0.1cm}}

\newtheorem{theor}{\noi Theorem}

\newtheorem{lemma}{\noi Lemma}

\newtheorem{corol}{\noi Corollary}

\newtheorem{prop}{\noi Proposition}
\newtheorem{defi}{\noi Definition}

\newtheorem{conj}{\noi Conjecture}

\theoremstyle{remark} 
\newtheorem{rem}{\noi Remark}
\newtheorem{rems}{\noi Remarks}
\newcommand{\ZZ}{{\mathbb Z}}
\newcommand{\QQ}{{\mathbb Q}} 
\newcommand{\NN}{{\mathbb N}} 
\newcommand{\RR}{{\mathbb R}}
\newcommand{\CC}{{\mathbb C}}
\DeclareMathOperator{\Hom}{Hom}

\newcommand{\G}{\mathbf G}

\newcommand{\proj}{\mathbf P}

\newcommand{\Ga}{\Gamma}
\newcommand{\La}{\Lambda}

\newcommand{\la}{\lambda}
\newcommand{\om}{\omega}

\newcommand{\R}{\mathbf R}
\newcommand{\OO}{\mathcal O}

\newcommand{\GL}{\mathbf{GL}}
\newcommand{\SL}{\mathbf{SL}}

\newcommand{\HHom}{\textnormal{Hom}}

\newcommand{\Ad}{\textnormal{Ad}}

\newcommand{\sh}{\textnormal{sh}}
\title{K\"ahler groups and duality}

\author{Bruno KLINGLER}

\begin{document}
\baselineskip 14pt
\setcounter{tocdepth}{1}

\maketitle
\tableofcontents
\section{Introduction and results}

\subsection{General setting}
The general setting for this paper is the study of topological
properties of compact K\"ahler manifolds, in particular complex smooth
projective varieties. The possible homotopy types for these spaces are
essentially unknown (c.f. \cite{si4}). An {\em a priori} simpler
question asks which finitely presented groups can be realized as fundamental groups of compact K\"ahler
manifolds (the so-called K\"ahler groups). 

It is well-known that any finitely presented group $\Ga$ can be realized as
the fundamental group of a $4$-dimensional compact real manifold, or
even of a symplectic $4$-manifold. A classical result of Serre
\cite{ser} shows
that any finite group can be realized as the fundamental group of a
smooth complex projective variety. On the other hand there are many known
obstructions for an infinite finitely presented group being K\"ahler (we refer to \cite{abckt}
for a panorama). Most of them come from Hodge theory (Abelian or not)
{\em in cohomological degree one}. 

As a prototype~: let $M$ be a compact K\"ahler
manifold with fundamental group $\Ga$. Classical Hodge theory shows
the existence of a weight $1$ pure Hodge structure on $H^1(\Ga, \QQ)=
H^1(M, \QQ)$. Thus $b^1(\Ga)$ has to be even. By considering finite
\'etale covers $b^1(\Ga')$ has to be even for any finite index
subgroup $\Ga'$ of $\Ga$. For example the free group on $n$-generator
$F_n$ is never K\"ahler, $n\geq 1$. 

The most interesting conjecture concerning infinite
K\"ahler groups, due to Carlson and Toledo and publicized by Kollar \cite{ko} and
Simpson, deals with {\em cohomology in degree $2$}~:
\begin{conj}[Carlson-Toledo] \label{conj1}
Let $\Ga$ be an infinite K{\"a}hler group. Then virtually $b^2(\Ga) >0$.
\end{conj}

\begin{rem} Recall that a group $\Ga$ has virtually some property $\mathcal{P}$ if
a finite index subgroup $\Ga' \subset \Ga$ has $\mathcal{P}$.
\end{rem}

Conjecture~\ref{conj1} means (c.f. appendix~\ref{explanation}) that there exists a compact K\"ahler manifold $M$
with $\pi_1(M) = \Ga'$ a finite index subgroup of $\Ga$ such that the
rational Hurewicz morphism $\pi_2(M) \otimes_\ZZ \QQ \lo H_2(M, \QQ)$
is not surjective. This statement is highly non-trivial~: there
exists compact K\"ahler manifolds $M$ (in fact smooth projective
complex varieties) whose $\pi_2(M)$ is very big, namely not 
finitely generated as a $\ZZ \Ga$-module \cite{DiPaSu}. 

Conjecture~\ref{conj1} trivially holds true for fundamental groups of
complex projective curves. By definition is is also satisfied by fundamental
groups of compact K\"ahler hyperbolic manifolds \cite{gromov}.
The strongest evidence in its favor is that it holds true if the
K\"ahler group $\Ga$ admits a finite dimensional complex representation $\rho$
with $H^1(\Ga, \rho) \not =0$ (c.f. theorem~\ref{prop2} appendix~\ref{sec2}). In
particular it is true in all the known examples with
very big $\pi_2(M)$ (they satisfy $b^1(\Ga)>0$).

\subsection{Results}
This paper is the first in a series of two studying
cohomological properties of K\"ahler groups. It is essentially
topological, using mainly duality in group cohomology and topological properties of Stein
spaces. The second one \cite{klingler} on the other hand is mainly geometric and
uses non-Abelian Hodge theory. One will also consult \cite{KlKoMau}
for a partial result towards conjecture~\ref{conj1}.

\sspace
Recall that a group $\Ga$ is an $r$-dimensional duality group (for
some positive integer $r$) if it satisfies a weak version of
Poincar\'e duality~: there exists a $\ZZ
\Gamma$-module $I$ such that for any $\ZZ \Gamma$-module $A$ there is
a natural isomorphism (i.e. induced by cap product with a fundamental
class)~:
$$ \forall i \in \NN, \quad H^i(\Ga, A) \simeq H_{r-i}(\Ga, I
\otimes_\ZZ A)\;\;,$$
where $I \otimes_\ZZ A$ denotes the tensor product over $\ZZ$ with
diagonal action.

Many groups of geometric origin are duality groups~: fundamental
groups of aspherical manifolds, arithmetic lattices, mapping class
groups .... 

\sspace
Our main result in this paper is the following~:

\begin{theor} \label{linear}
Let $\Ga$ be an infinite linear $r$-dimensional duality group, $r \geq
6$. If $\Ga$ is a K\"ahler group then virtually
$$b^2(\Ga) + b^4(\Ga) >0 \;\;.$$
\end{theor}

\begin{rem}
Requiring the group $\Ga$ to be linear is quite restrictive as we know that there
exist many fundamental groups of smooth complex projective
varieties which are not linear, not even residually finite
\cite{to}. This assumption is relaxed in \cite{klingler} under some
unboundedness condition.
\end{rem}

Theorem~\ref{linear} excludes many groups
from being K\"ahler. A striking example is the following~:

\begin{theor} \label{p-adic}
Let $G_v$ be the group of $F_v$-points of an algebraic
group $\G_v$, with reductive neutral component, over a non-Archimedean
local field $F_v$ of characteristic $0$. 
Suppose that $\text{rank}_{F_{v}} \G_v \geq 6$.

Then a cocompact lattice $\Ga \subset G_v$ is not K\"ahler.
\end{theor}

\begin{rem}
Notice that the conclusion of theorem~\ref{p-adic} is predicted
by Simpson's integrality conjecture under the weaker assumption
$\text{rank}_{F_{v}} \G_v \geq 2$ (c.f. section~\ref{padic}).
\end{rem}

The proof of theorem~\ref{linear} relies on the following
ideas. First the group theoretical properties of the fundamental group
of a space $M$ are intimately linked to the geometric properties of
the universal cover~$\ti{M}$ of $M$. In the case of K\"ahler manifolds
one has the following famous conjecture~:

\begin{conj}[Shafarevich] \label{shafa}
Let $M$ be a connected compact K\"ahler manifold. Then its universal cover
$\ti{M}$ is holomorphically convex.
\end{conj}

\begin{rem}
We recall that a complex space $E$ is holomorphically convex if
for any sequence $(x_n)_{n \in \NN}$ of points in $E$
there exists a sequence of holomorphic functions $(f_n \in \mathcal{O}_E(E))_{n
  \in \NN}$ such that the complex sequence $(f_n(x_n))_{n \in \NN}$ is
unbounded. Equivalently $E$ admits a proper map $E \lo S$ with $S$
Stein, $S$ is called the Cartan-Remmert reduction of $E$.
\end{rem}

From this point of view it is natural to first study the fundamental
groups of compact K\"ahler manifolds $M$ whose universal cover
$\ti{M}$ is Stein. If one moreover assumes that $M$ is a smooth
projective complex variety one can assume that $M$ is a smooth
projective complex surface using Lefschetz's hyperplane theorem. Using the
topological restrictions on $\ti{M}$ coming from the Stein condition, the
Leray spectral sequence comparing the cohomologies of $M$ and $\Ga$,
and the fact that $\Ga$ and $M$ both satisfy duality but in different degrees, one shows the following~:

\begin{theor} \label{stein}
Let $\Ga$ be an infinite oriented $r$-dimensional duality group, $r \geq 6$.
Suppose $\Ga = \pi_1(M)$ where $M$ is a connected smooth projective complex
variety with Stein universal cover $\ti{M}$. Then
$$b^1(\Ga) + b^4(\Ga) >0\;\;.$$
\end{theor}

\begin{rem}
Notice that in theorem~\ref{stein} we are not asking for $\Ga$ to be
linear. Notice also that in the conclusion of theorem~\ref{stein} we
do not need to pass to a finite index subgroup.
\end{rem}

To deduce theorem~\ref{linear} from theorem~\ref{stein} one uses
non-Abelian Hodge theory as in the proof of the Shafarevich conjecture
for linear fundamental groups given in \cite{katram} and
\cite{eyssi}. One is essentially reduced to proving a statement
similar to theorem~\ref{stein} with $M$ only normal. Although I can't
prove theorem~\ref{stein} in the normal case the extra 
informations in the case at hand are enough to conclude the proof of theorem~\ref{linear}.

\subsection{Acknowledgements}
It is a pleasure to thank Pierre Vogel for topological
discussions, P.Eyssidieux,  J.Maubon, C. Simpson and D.Toledo for their interest
in this work. I especially thanks V. Koziarz and B.Claudon for bringing \cite{zuo1}
to my attention, thus extending my results from fundamental groups of
smooth projective complex varieties to K\"ahler groups.

\subsection{Notations}
If $M$ is a topological space we abbreviate by $H_\bullet(M)$
(resp. $H^\bullet(M)$) its integral homology (resp. cohomology)
$H_\bullet(M, \ZZ)$ (resp. $H^\bullet(M, \ZZ)$).

\section{Duality and finiteness}

The goal of this section is to prove the following two lemmas
(probably well-known to the algebraic topologist). 

The first lemma implies that the duality requirement for $\Ga$ excludes the worse possible situation
for conjecture~\ref{conj1}~:
\begin{lemma} \label{finite}
Let $\Ga$ be a duality group. Suppose $\Ga =\pi_1(M)$ for $M$ a
connected compact manifold. Then $\pi_2(M)$ is a finitely generated (left)
$\ZZ \Ga$-module.
\end{lemma}

The second lemma implies that given a closed oriented connected $4$-manifold with fundamental
group $\Ga$, the $\ZZ \Gamma$-module $\pi_2(M)$ satisfies a strong duality property
as soon as $\Ga$ is an $r$-dimensional duality group, $r >3$~:

\begin{lemma} \label{pi2dual}
Let $M$ be a closed oriented connected $4$-manifold with fundamental
group $\Ga$. If $H^2(\Ga,\Lambda_\ZZ)
= H^3(\Ga, \Lambda_\ZZ) =0$ then one has a natural isomorphism of (left) $\Lambda_\ZZ$-modules
$$ \pi_2(M) \simeq \overline{\HHom_{\La_{\ZZ}}(\pi_2(M),
  \Lambda_\ZZ}) 
$$
(where $\Lambda_\ZZ$ denotes the group ring $\ZZ \Gamma$ and if $R$ is a
right $\Lambda_\ZZ$-module the notation $\overline{R}$ denotes the
associated left $\Lambda_\ZZ$-module).
\end{lemma}

\subsection{Finiteness conditions}
We first recall some basic definitions before
proving lemmas~\ref{finite} and \ref{pi2dual}.

\begin{defi}
A group $G$ is said to be of type $F_n$, $0 \leq n < + \infty$, if
there is a $K(G,1)$ with a finite $n$-skeleton, i.e. with only
finitely many cells in dimensions $\leq n$. We say that $G$ is of type
$F_{\infty}$ if there is a $K(G,1)$ with all its skeleta finite. We
say that $G$ is of type $F$ if there is a finite $K(G,1)$.
\end{defi}

Obviously any group if of type $F_0$, a group $G$  is of type $F_1$ if
and only if it is finitely generated and of type $F_2$ if and only if
it is finitely presented.

\begin{defi}
A group $G$ is said to be of type $FP_n$, $0 \leq n < + \infty$, if
there is a projective resolution $P_\bullet$ of $\ZZ$ over $\ZZ G$
such that $P_i$ is finitely generated for $i \leq n$. We say that $G$
is of type $FP_{\infty}$ if there is a projective resolution
$P_\bullet$ of $\ZZ$ over $\ZZ G$ with $P_i$ finitely generated for
all $i$. We say that $G$ is of type $FP$ if there is a finite
projective resolution $P_\bullet$ of $\ZZ$ over $\ZZ G$.
\end{defi}

Obviously a group $G$ of type $F_n$ is of type $FP_n$ as the cellular
chain complex $C_\bullet(\widetilde{K(G,1)})$ (with $\ZZ$-coefficients) provides a free (thus
  projective) resolution of $\ZZ$ as a $\ZZ G$-module up to degree $n$. One can show that
$G$ is of type $F_1$ if and only if it is of type $FP_1$ and that $G$
is of type $FP_2$ if and only if $G = \tilde{G}/N$ where $\tilde{G}$
is of type $F_2$ and $N$ is a perfect normal subgroup of
$\tilde{G}$. In this case $G$ is of type $F_2$ if and only if $N$ is
finitely generated as a normal subgroup. Bestvina and Brady exhibited
examples where this is not the case \cite{bestbra}. For higher $n$
however a group $G$ is of type $F_n$, $3 \leq n \leq \infty$, if and
only if $G$ is finitely presented and of type $FP_n$.

\subsection{Duality groups}
If $\Ga$ is the fundamental group of an aspherical closed $n$-manifold
$M$ then if satisfies Poincar\'e duality~:
$$H^i(\Ga, A) \simeq H_{n-i}(\Ga, I \otimes_\ZZ A)$$
where $I$ denotes the orientation  module of $M$ and $A$ can be any $\ZZ
\Gamma$-module. This leads Bieri \cite{Bie} to the following~:

\begin{defi} \label{poincare}
A group $\Ga$ is a Poincar\'e duality group of dimension $n$ it there
is a $\ZZ \Gamma$-module $I$ (called the dualizing module of $\Ga$), which is isomorphic to $\ZZ$ as a
$\ZZ$-module, and a homology class $\mu \in H_n(\Ga, I)$ so that for
any $\ZZ \Gamma$-module $A$, cap-product with $\mu$ defines an
isomorphism~:
$$H^i(\Ga, A) \simeq H_{n-i}(\Ga, I \otimes_\ZZ A)\;\;.$$
\end{defi}

Since the universal covering $\ti{M}^n$ of an aspherical manifold $M^n$
is contractible it follows from Poincar\'e duality in the non-compact
case that the cohomology with compact support of $\ti{M}^n$ is the
same as that of $\RR^n$ i.e.~:
$$H^i_c(\ti{M}^n) = \begin{cases}
0, \quad \text{for} \; i \not=n, \\
\ZZ, \quad \text{for} \; i =n.
\end{cases}
$$
On the other hand $\Ga$ acts freely properly and cocompactly on the
acyclic space $\ti{M}^n$ thus $H^i(\Ga, \ZZ \Gamma) \simeq
H^i_c(\ti{M}^n)$. Johnson and Wall \cite{JoWa} proposed the following
definition, which was shown to be equivalent to the previous one in
\cite{BieEck}~:
\begin{defi} \label{poincarebis}
A group $\Ga$ is a Poincar\'e duality group of dimension $n$ if $\Ga$
is of type $FP$ and $$H^i(\Ga, \ZZ \Ga) = \begin{cases}
0, \quad \text{for} \; i \not=n, \\
\ZZ, \quad \text{for} \; i =n.
\end{cases}
$$
\end{defi}
\noi
The dualizing module is $I = H^n(\Ga, \ZZ \Ga)$. 
The main content of the equivalence of definitions~\ref{poincare} and \ref{poincarebis} is that
definition~\ref{poincare} forces $\Ga$ to be of type $FP$.

There are many interesting groups which satisfy
definition~\ref{poincare} except for the requirement that $I$ is
isomorphic to $\ZZ$ (as a $\ZZ$-module). This leads to the following
theorem-definition \cite{BieEck}~:
\begin{theor}[Bieri-Eckmann] \label{dual}
A group $\Ga$ is a duality group of dimension $n$ if it satisfies one
the following equivalent conditions~:
\begin{itemize}
\item There exists a $\ZZ \Gamma$-module $I$ such that for any $\ZZ
  \Gamma$-module $A$ there is an isomorphism induced by cap-product
  with a fundamental class~: $H^i(\Ga, A) \simeq H_{n-i}(\Ga, I
  \otimes_\ZZ A)$.
\item $\Ga$ is of type $FP$ and 
$$H^i(\Ga, \ZZ \Ga) = \begin{cases}
0, \quad \text{for} \; i \not=n, \\
I, \quad \text{for} \; i =n.
\end{cases}
$$
\end{itemize}
\end{theor}

\subsection{A finiteness lemma}
Lemma~\ref{finite} is an immediate corollary of the following more precise~:
\begin{lemma}  \label{finite1}
Let $M$ be a connected $CW$-complex with fundamental group
$\Gamma$. Let $\La_\ZZ= \ZZ \Ga$ be the group ring of $\Gamma$. 

Suppose $M$ has a finite $3$-skeleton.

Then $\pi_2(M)$ is finitely generated as a $\La_\ZZ$-module if and only if
$\pi_1(M)$ is of type $FP_3$.
\end{lemma}

\begin{proof}
Suppose first that $\pi_2(M)$ is finitely generated as a
$\La_\ZZ$-module. Thus one can construct the $3$-squeletton of a $K(\Ga, 1)$
by adding finitely many $3$-cells to the $3$-skeleton $M^{(3)}$ of
$M$. As $M^{(3)}$ is finite the group $\Ga$ is of type $F_3$ thus also $FP_3$.

Conversely suppose $\Ga$ is of type $FP_3$. Thus there is a resolution
$(P_\bullet,d)$ of the $\La_\ZZ$-module $\ZZ$ with $P_i$ a projective finitely
generated $\La_\ZZ$-module, $i \leq 3$. Let $B_2 = d(P_3)$ the group of
$2$-boundaries. As $P_3$ is $\La_\ZZ$-finitely generated, $B_2$ too and
one has the following exact sequence of finitely generated $\La_\ZZ$-modules~:
\begin{equation} \label{eq1}
0 \lo B_2 \lo P_2 \lo P_1 \lo P_0\lo \ZZ \lo 0\;\;.
\end{equation}
One the other hand the cellular chain complex of the
universal cover $\tilde{M}$ gives the following exact sequence of $\La_\ZZ$-modules~:
\begin{equation} \label{eq2}
0 \lo Z_{2}(\tilde{M}) \lo C_2(\ti{M}) \lo C_1(\ti{M}) \lo
C_0(\ti{M}) \lo \ZZ \lo 0 \;\; ,
\end{equation}
where $Z_2(\ti{M}) = \ker (d_2: C_2(\ti{M}) \lo C_1(\ti{M}))$.
As $M$ has a finite $3$-skeleton the $C_i(\ti{M})$'s, $0 \leq i \leq 3$,
are finitely generated $\La_\ZZ$-modules.

Recall now the following classical
\begin{lemma}[Schanuel]
Let $R$ be a ring. If 
$$0 \lo B_n \lo P_{n-1} \lo \cdots \lo P_1 \lo
P_0 \lo A \lo 0$$
and 
$$0 \lo B'_n \lo P'_{n-1} \lo \cdots \lo P'_1 \lo P'_0 \lo
A \lo 0$$
are exact sequences of $R$-modules with $P_i, P'_i$ projective
$R$-modules $(0 \leq i \leq n-1)$ then~:
$$ B_n \oplus P'_{n-1} \oplus P_{n-2} \oplus \cdots \simeq B'_n \oplus
P_{n-1} \oplus P'_{n-2} \oplus \cdots \;\;.$$
\end{lemma}

Applying Schanuel's lemma to (\ref{eq1}) and (\ref{eq2}) one obtains~:
$$B_2 \oplus C_2(\ti{M}) \oplus P_1 \oplus \ZZ \simeq Z_2(\ti{M})
\oplus P_2 \oplus C_1(\ti{M}) \oplus P_0 \oplus \ZZ\;\;.$$
This implies that $Z_2(\ti{M})$ is a finitely generated $\La_\ZZ$-module,
thus also its quotient $\pi_2(M) = H_2(\ti{M}) = Z_2(\ti{M}) /
B_2(\ti{M})$.
\end{proof}

\subsection{Topology of closed $4$-manifolds and proof of lemma~\ref{pi2dual}}

Let $M$ be a closed oriented connected $4$-manifold with fundamental group
$\Ga$ and universal cover $\tilde{M}$. Denote by $\La_\ZZ = \ZZ \Ga$ the
group ring of $\Gamma$. An interesting invariant is the 
equivariant intersection form~:
$$ s_M : \pi_2(M) \otimes_\ZZ \pi_2(M) \lo \La_\ZZ \;\;,$$
generalizing the well-known intersection product 
$$H_2(M, \ZZ) \otimes_\ZZ H_2(M, \ZZ) \lo \ZZ$$ 
and defined as follows.

By Poincar\'e duality there is a canonical isomorphism $H^i_c(\tilde{M}, \ZZ) \simeq H_{4-i}(\tilde{M}, \ZZ)$
between cohomology with compact support and homology. In particular~:
\begin{itemize}
\item[a)]  this induces an isomorphism
$$\varepsilon : H^4_c(\tilde{M} , \ZZ) \simeq  H_0(\tilde{M}, \ZZ)
=\ZZ\;\;.$$
\item[b)] 
Let $h : \pi_2(M) \stackrel{\sim}{\lo} H_2(\tilde{M}, \ZZ)$ be the
Hurewicz isomorphism. Using Poincar\'e duality we obtain a 
canonical isomorphism $$ \xymatrix@1{\pi_2(M) \ar[r]^\sim_{\phi}
  &H^2_c(\tilde{M}, \ZZ)} \;\;.$$
\end{itemize}
The cup-product on $H^2_c(\tilde{M}, \ZZ)$ is a $\Ga$-invariant linear
form with values in $\ZZ$. We thus define
$$s_M (x, y) = \sum_{g \in \Gamma} \varepsilon(\phi(g^{-1}x) \cup
\phi(y))\cdot g  \in \La_\ZZ \;\;.$$
Notice that this pairing is $\La_\ZZ$-hermitian in the sense that for all
$\la \in \La_\ZZ$ we have
$$s_M (\la \cdot x, y) = \lambda \cdot s_M(x,y) \quad \text{and}
\quad s_M(y,x) = \overline{s_M(x,y)} \;\;,$$
where $\La_\ZZ$ acts on itself by left translation and the involution $\la
\mapsto \overline{\la}$ on $\La_\ZZ$ is given by $\overline{g} = g^{-1}$
for $g \in \Ga$.

\sspace
This pairing is controlled as follows. The Leray-Serre cohomological
spectral sequence with coefficients in $\La_\ZZ$ for the classifying map
$M \lo K(\Ga, 1)$, whose homotopy fiber is $\tilde{M}$, yields the
following short exact sequence of left $\La_\ZZ$-modules~:
$$ 0 \lo H^2(\Ga, \La_\ZZ) \lo H^2(M, \La_\ZZ) \lo \overline{\HHom_{\La_\ZZ}(H_2(M, \La_\ZZ),
\La_\ZZ)} \lo H^3(\Ga, \La_\ZZ) \lo 0\;\;.$$

By Poincar\'e duality for $M$ and Hurewicz theorem we have an isomorphism of
left $\La_\ZZ$-modules~:
$$H^2(M, \La_\ZZ) \simeq H_2(M, \La_\ZZ) \simeq H_2(\tilde{M}, \ZZ) \simeq
\pi_2(M) \;\;.$$
Finally we get the following exact sequence of left $\La_\ZZ$-modules~:
\begin{equation} \label{eq3}
0 \lo H^2(\Ga, \La_\ZZ) \lo \pi_2(M) \stackrel{\tau_M}{\lo}
\overline{\HHom_{\La_\ZZ}(\pi_2(M), \La_\ZZ)} \lo H^3(\Ga, \La_\ZZ) \lo 0\;\;.
\end{equation}
One easily checks that the map $\tau_M$ is the
natural map associated to the pairing $s_M$. In particular the kernel
of $s_M$ is nothing else than $H^2(\Ga, \La_\ZZ)$, its cokernel is
$H^3(\Ga, \La_\ZZ)$. 

From this short discussion we immediately get lemma~\ref{pi2dual}.

%%%%%%%%%%%%%%%%%%%%%%%%%%%%%%%%%%%%%%%%%

\section{Cohomological properties of fundamental groups of smooth
  projective surfaces
  whose universal cover is Stein} 

\begin{prop} \label{leraycoh}
Let $M$ be a connected $4$-dimensional $CW$-complex with universal
cover $\ti{M}$ and fundamental group $\Ga$. Let $R$ be any (left)
$\Ga$-module. If $H^3(\ti{M}) = H^4(\ti{M})=0$ then
\begin{itemize}
\item[(a)] the following exact sequence of $\Ga$-modules holds~:
\begin{equation*} \label{excoh}
\begin{split}
0 \lo H^2(\Ga, R) \lo H^2(M, R) \lo (H^2(\ti{M}) \otimes_\ZZ R)^\Gamma
\lo H^3(\Ga, R) \lo H^3(M, R) \lo \\
H^1(\Ga, H^2(\ti{M}) \otimes_\ZZ R)
\lo H^4(\Ga, R) \lo H^4(M, R) \lo H^2(\Ga, H^2(\ti{M}) \otimes_\ZZ R)
\lo H^5(\Ga, R) \lo 0\;\;.
\end{split}
\end{equation*}
\item[(b)] $ \forall\; i \geq 3, \quad H^i(\Ga, H^2(\ti{M}) \otimes_\ZZ R)
  \simeq H^{i+3}(\Ga, R) \;\;.$
\end{itemize}
\end{prop}

\begin{proof}
Let us write the Leray-Serre spectral sequence associated to the
classifying map $f: M \lo K(\Ga,1)$ (whose homotopy fiber is
$\tilde{M}$) for the local coefficients~$R$~:
$$ 
\xymatrix@1{H^p(\Ga, H^q(\tilde{M}) \otimes_\ZZ R) \ar@{=>}[r] & H^{p+q}(M,
  R)} \;\;.$$
As $\tilde{M}$ is simply connected, the $E_2$-page of this spectral
sequence has only two non-zero lines~:
$$E_2^{p, 0} = H^p(\Ga, R) \qquad \text{and} \qquad E_2^{p,2} =
H^p(\Ga, H^2(\tilde{M}) \otimes_\ZZ R)\;\;.$$

It implies that $E_3=E_2$ and then $E_4 = E_5 = \cdots = E_{\infty}$.
Thus~:
\begin{gather*}
E_{\infty}^{0, 0} = R^\Gamma  \\
E_{\infty}^{0,1} = 0, \quad  E_{\infty}^{1,0} =
H^1(\Ga, R), \\
E_{\infty}^{0,2} = \ker ((H^2(\ti{M}) \otimes_\ZZ R)^\Ga \stackrel{d_3}{\lo}
H^3(\Ga, R)) , \quad E_\infty^{1,1} = 0, \quad E_\infty^{2, 0}= H^2(\Ga, R),
\end{gather*}
and for all $i, j$, $i+j \geq 3$, 
$$
E_{\infty}^{i, j} = \begin{cases}
\ker(H^i(\Ga, H^2(\ti{M}) \otimes_\ZZ R) \stackrel{d_3}{\lo}
H^{i+3}(\Ga, R)) & \text{if} \; j=2, \\
H^i(\Ga, R)/d_3 (H^{i-3}(\Ga, H^2(\ti{M}) \otimes_\ZZ R))& \text{if} \; j = 0, \\
0 & \text{otherwise.}
\end{cases}
$$

From this one deduces the following short exact sequences~:
$$ 0 \lo H^2(\Ga, R) \lo H^2(M, R) \lo \ker((H^2(\ti{M}) \otimes_\ZZ R)^\Ga
\stackrel{d_3}{\lo} H^3(\Ga, R)) \lo 0\;\;,$$
and for $i \geq 3$~:
$$0 \lo \frac{H^i(\Ga, R)}{d_3(H^{i-3}(\Ga,
H^2(\ti{M}) \otimes_\ZZ R))}\lo H^i(M, \CC) \lo \ker(H^{i-2}(\Ga,
H^2(\ti{M}) \otimes_\ZZ R) \stackrel{d_3}{\lo} H^{i+1}(\Ga, R)) \lo 0\;\;.$$
The proposition then follows from the consideration of these short
exact sequences for $i \leq 5$ and the fact that $H^i(M, \cdot) = 0$
for $i \geq 5$.
\end{proof}

Studying similarly the Leray-Serre homology spectral sequence
one obtains the dual statements~:
\begin{prop} \label{lerayho}
Let $M$ be a connected $4$-dimensional $CW$-complex with universal
cover $\ti{M}$ and fundamental group $\Ga$. Let $R$ be any (left)
$\Ga$-module. If $H_3(\ti{M})=H_4(\ti{M})=0$ then~:
\begin{itemize}
\item[(a)] the following exact sequence of $\Ga$-modules holds~:
\begin{equation*} \label{exho}
\begin{split}
0 \llo H_2(\Ga, R) \llo H_2(M, R) \llo (H_2(\ti{M}) \otimes_\ZZ R)_\Gamma
\llo H_3(\Ga, R) \llo H_3(M, R) \llo \\
H_1(\Ga, H_2(\ti{M}) \otimes_\ZZ R)
\llo H_4(\Ga, R) \llo H_4(M, R) \llo H_2(\Ga, H_2(\ti{M}) \otimes_\ZZ R)
\llo H_5(\Ga, R) \llo 0\;\;.
\end{split}
\end{equation*}
\item[(b)] $ \forall\; i \geq 3, \quad H_i(\Ga, H_2(\ti{M}) \otimes_\ZZ R)
  \simeq H_{i+3}(\Ga, R) \;\;.$
\end{itemize}
\end{prop}

\begin{corol} \label{stein0}
Let $M$ be a connected complex surface (not necessarily smooth) with
Stein universal cover $\ti{M}$ and fundamental group $\Ga$. Let $R$ be any (left)
$\Ga$-module. Then the conclusions of proposition~\ref{leraycoh} and
proposition~\ref{lerayho} hold true.
\end{corol}

\begin{proof}
The fact that a Stein surface $\ti{M}$ (not necessarily smooth)
satisfies $H^i(\ti{M})=0$, $i \geq 3$, is a classical result of
Narasimhan \cite{nar}. In fact such an $\ti{M}$ is homotopy equivalent
to a $2$-dimensional $CW$-complex \cite{GorMac}.
\end{proof}

\section{On fundamental groups of smooth projective
  varieties whose universal cover is Stein~: proof of theorem~\ref{stein}}

\subsection{A vanishing lemma}
\begin{lemma} \label{invariant}
Let $\Ga$ be an infinite group of type $FP_3$ satisfying $$H^2(\Ga, \La)
= H^3(\Ga, \La) =0\;\;,$$
where $\La = \QQ \Gamma$.
Suppose $\Ga$ is the
fundamental group of a connected $CW$-complex $M$. Then 
$$(H^2(\ti{M})_\QQ \otimes_\QQ \pi_2(M)_\QQ)^\Ga = 0$$
(for the diagonal action of $\Ga$ on $H^2(\ti{M}) \otimes_\ZZ \pi_2(M)$).
\end{lemma}

\begin{proof}
By lemma~\ref{finite1} $\pi_2(M)$ is a finitely generated (left)
$\ZZ \Ga$-module. Thus there exists a positive integer $i$ and a
surjective morphism of (left) $\Lambda$-modules~:
$$ \Lambda^i \twoheadrightarrow \pi_2(M)_\QQ \;\;.$$

This implies that one has an injective morphism of (left) $\Lambda$-modules~:
$$\overline{\HHom_{\La}(\pi_2(M)_\QQ, \La)} \hookrightarrow \Lambda^i\;\;.$$

By lemma~\ref{pi2dual}~:
$$ \pi_2(M) \simeq \overline{\HHom_{\La_{\ZZ}}(\pi_2(M),
  \Lambda_\ZZ}) \;\;.
$$
 
Finally $H^2(\ti{M})_\QQ \otimes_\QQ \pi_2(M)_\QQ \subset
H^2(\ti{M})_\QQ \otimes_\QQ \Lambda^{i}$ and one is reduced to
proving~: for any non-trivial (left) $\Lambda$-module $R$ one has
$$(R \otimes_\QQ \Lambda)^\Ga = 0\;\;.$$
This follows from the fact that the infinite group $\Ga$ has only
infinite orbits in its action by left translation on $\La = \QQ \Ga$.

\end{proof}
%%%%%%%%%%%%%%%%%%%%%%%%%%%%
\subsection{The topological theorem}

\begin{theor} \label{top}
Let $\Ga$ be an infinite oriented $r$-dimensional duality group, $r \geq 6$.
Suppose $\Ga = \pi_1(M)$ where $M$ is a connected oriented compact
$\mathcal{C}^\infty$ $4$-manifold whose universal cover $\ti{M}$
has the homotopy type of a $2$-dimensional $CW$-complex. Then
$$b^1(\Ga) + b^4(\Ga) >0\;\;.$$
\end{theor}

\begin{proof}
Assume by contradiction $b^1(\Ga) + b^4(\Ga) =0$.
As $b^1(M)= b^1(\Ga)$ and $b^3(M) = b^1(M)$ by Poincar\'e duality for
$M$, one obtains $b^3(M)=0$.

The homological exact sequence~(proposition~\ref{lerayho}, (a))
gives the exact sequence~:
$$H_4(\Ga)_\QQ \lo H_1(\Ga, \pi_2(M)_\QQ) \lo H_3(M)_\QQ $$
from which we deduce~:
\begin{equation} \label{e2}
H_1(\Ga, \pi_2(M)_\QQ)=0\;\;.
\end{equation}

By Poincar\'e duality for $M$ one deduces from equation~(\ref{e2})~:
\begin{equation} \label{e3}
H^3(M, \pi_2(M)_\QQ) \simeq H_1(M, \pi_2(M)_\QQ)=0\;\;.
\end{equation}

Let us now consider the cohomological exact sequence~(proposition~\ref{leraycoh},
(a)) with coefficients $\pi_2(M)_\QQ$~:
\begin{equation} \label{e4}
(H^2(\ti{M})_\QQ \otimes_\QQ \pi_2(M)_\QQ)^\Ga \lo H^3(\Ga,
\pi_2(M)_\QQ) \lo H^3(M, \pi_2(M)_\QQ) \;\;.
\end{equation}
As $\Ga$ is an $r$-dimensional duality group with $r \geq 3$~:
$$H^2(\Ga, \La) = H^3(\Ga, \La) =0\;\;.$$
By lemma~\ref{invariant} $(H^2(\ti{M})_\QQ \otimes_\QQ \pi_2(M)_\QQ)^\Ga=0$. On the other
hand by equation~(\ref{e3}) one gets $H^3(M, \pi_2(M)_\QQ) =0$. We
deduce from equation~(\ref{e4})~:
\begin{equation} \label{e5}
H^3(\Ga, \pi_2(M)_\QQ) =0\;\;.
\end{equation}

However by duality for $\Ga$~:
$$H^3(\Ga, \pi_2(M)_\QQ) \simeq H_{r-3}(\Ga, I \otimes_\ZZ\pi_2(M)_\QQ) \;\;,$$
and by the homological isomorphism~(proposition~\ref{lerayho}, (b)) one obtains as
$r \geq 6$~:
$$H_{r-3}(\Ga, I \otimes\pi_2(M)_\QQ) \simeq H_r(\Ga, I) \simeq
H^0(\Ga, \QQ) \simeq \QQ\;\;.$$
Finally $H^3(\Ga, \pi_2(M)_\QQ) \simeq \QQ$, contradiction to
equation~(\ref{e5}).
\end{proof}
%%%%%%%
\subsection{Proof of theorem~\ref{stein}}

Let $\Ga$ be a group as in theorem~\ref{stein}. Thus $\Ga= \pi_1(M)$,
$M$ complex smooth projective variety of (complex) dimension $n$ with
$\ti{M}$ Stein. If $n=1$ the variety $M$ is a smooth projective curve, of genus $\geq
1$ as $\Ga$ is infinite. In particular $M$ is a $K(\Ga, 1)$ and $b^1(\Ga) = b^1(M) >0$ thus
the conclusion of theorem~\ref{linear} is valid in this case. Therefore we
can assume $n\geq 2$. 

Fix $\OO_M(1)$ an ample line bundle on $M$. By Lefschetz hyperplane
theorem any smooth surface $N$ complete intersection in $M$ of hyperplane
sections corresponding to $\OO_M(1)$ still has fundamental group
$\Ga$. Moreover the universal cover $\tilde{N}$ of such an
$N$ is an irreducible component of the preimage of $N$ in
$\tilde{M}$. Thus $\ti{N}$ is a closed analytic submanifold of a Stein
manifold, hence a Stein manifold. Replacing $M$ by $N$  we can assume $n=2$.

Finally $M$ satisfies the hypotheses of theorem~\ref{top} and we conclude.

\diam

%%%%%%%%%%%%%%%%%%%%%%%%%%%%
%%%%%%%%%%%%%%%%%%%%%%%%%%%%
\section{Proof of theorem~\ref{linear}}

%%%%%%%%%%%%%%%%%
\subsection{Shafarevich maps}
 In \cite[theor.2]{eyssi}  Eyssidieux shows the following
\begin{theor}[Eyssidieux] \label{eyssi1}
Let $X$ be a (connected) smooth complex
  projective variety. Let $\rho : \pi_1(X) \lo GL(n, \CC)$
  be a semi-simple representation. 

There exists a diagram of analytic morphisms
$$ 
\xymatrix{
\tilde{X}/ \ker \rho \ar[d] \ar[r]^{\widetilde{\sh_\rho}} &
\widetilde{S_\rho}(X) \ar[d]_{\pi} \\
X \ar[r]_{\sh_{\rho}} & \sh_\rho(X)}
$$
where~:
\begin{itemize}
\item $\sh_\rho(X)$ is a normal projective variety, the morphism
  $\sh_\rho: X \lo \sh_\rho(X)$ has connected fibers and for any
morphism $Z \lo X$, with $Z$ a smooth connected projective variety, the image
$\sh_\rho(Z)$ is a point if and only if $\rho(\pi_1(Z))$ is finite.
\item $\widetilde{S_\rho}(X)$ is a normal analytic variety without compact
  positive dimensional analytic subspaces with a proper discontinuous
  action of $\rho(\Ga)$, the morphism $\tilde{X}/ \ker \rho \lo
  \widetilde{S_\rho}$ is $\Ga$-equivariant, proper and satisfies
  $\sh_\rho(X) = \widetilde{S_\rho}(X)/\rho(\Gamma)$.
\end{itemize}
\end{theor}

In the special case where $\rho$ is rigid Eyssidieux proves the
holomorphic convexity of the cover~$\tilde{X}/\ker \rho$
(c.f. \cite[theor.3]{eyssi} and \cite[section 4]{eyssi} for the proof)~:

\begin{theor}[Eyssidieux] \label{eyssi2}
With the notations and assumptions of theorem~\ref{eyssi1} suppose
moreover that $\rho$ is rigid. Then $\widetilde{S}_\rho(X)$ is Stein
(and thus $\ti{X}/\ker \rho$ is holomorphically convex).
\end{theor}

\begin{rems}
\begin{itemize}
\item[(a)] In addition to \cite[theor.2]{eyssi} the reader might want
  to look at \cite[theor.2.1.7]{eyssi} and \cite[prop.2.2.20]{eyssi}
  which deal with the case of a rigid $\rho$ (the only case of
  interest for us).
\item [(b)]
In the case where $X$ is a surface,
theorem~\ref{eyssi1} and theorem~\ref{eyssi2} in this context are proven in
\cite{katram}. However we refer to \cite{eyssi} for a more
detailed construction.
\item[(c)] Of course the results of \cite{eyssi} are more general. In
  particular the holomorphic convexity of $\ti{X}/\ker \rho$
  generalizes to the case where one replaces $\rho$ by a quasi compact
  absolutely constructible set $R$ of conjugacy classes of semisimple
  representations in the sense of Simpson and $\tilde{X}/\ker \rho$ by
  $\tilde{X}_R$ the cover of $X$ corresponding to the intersection of
  their kernels.
\end{itemize}
\end{rems}

We will need one more fact from~\cite{eyssi}. Suppose $\rho$ is a
rigid representation. Thus one can assume that $\rho$ takes values in $\GL(n, K)$, $K$ number
field. Let $\G$ be the natural $K$-form of the reductive group
Zariski closure of $\rho(\Ga)$ in $\GL(n, \CC)$. As $\Ga$ is of finite
type one has $\rho(\Ga) \subset \G((\OO_{K})_S)$ where $\OO_K$ denotes
the ring of integers of $K$ and $S$ is a finite set of places of
$K$. We assume $S$ minimal i.e. for any $v \in S$ the representation
$\rho_v : \Ga \lo \G(K_v)$ does not have bounded image.
Let $S_{\text{ar}} \subset S$ the set of archimedean places an $S_f =
S \setminus S_{\text{ar}}$ the set of finite places. As in
\cite[p.524]{eyssi} let 
$$h : \ti{X}/\ker \rho \lo R_S:=\prod_{v \in S_{f}} \Delta_v \times
\prod_{v \in S_{\text{ar}}} R_v$$ 
be the natural $\Ga$-equivariant harmonic map to the product of the Bruhat-Tits buildings $\Delta_v$ associated to
the $p$-adic groups $\G_{K_{v}}$, $v \in S_f$, and of the symmetric
spaces $R_v$ associated to the real Lie group $\G_{K_{v}}$, $v \in
S_{\text{ar}}$. Then by construction of $\widetilde{S}_{\rho}(X)$ one
has (c.f. \cite[p.524]{eyssi})~:
\begin{lemma}[Eyssidieux] \label{leEyssi}
The harmonic map $h$ factorizes (equivariantly) through $\widetilde{S}_{\rho}(X)$~:
$$ 
\xymatrix{
\ti{X}/\ker \rho \ar[r]^h \ar[d]_{\widetilde{\sh_{\rho}}}  &R_S \\
\widetilde{S}_{\rho}(X) \ar@{-->}[ru] &}\;\;.$$
\end{lemma}

\subsubsection{Proof of theorem~\ref{linear}}

Let $\Ga$ be a group as in theorem~\ref{linear}. Thus $\Ga= \pi_1(X)$,
$X$ compact K\"ahler manifold. Let $\rho: \Ga \hookrightarrow GL(n,
\CC)$ be a linear embedding. 

Suppose by contradiction $b^2(\Ga) + b^4(\Ga) =0$. As explained in the appendix~\ref{sec2}
our assumption $b^2(\Ga) =0$ implies that $\Ga$ is {\em
  schematically rigid}, meaning that for any $n \in \NN$
and any representation $\rho: \Ga \lo \GL(n, \CC)$ one has $H^1(\Ga,
\Ad \rho)=0$. Considering the trivial representation we obtain in
particular that 
\begin{equation} \label{e7}
b^1(\Ga)=0\;\;.
\end{equation}
Schematic rigidity implies that $\Ga$ is {\em reductive}, meaning that all its linear finite-dimensional
complex representations are semi-simple (c.f. \cite[section
3.5]{KlKoMau}). In particular the Zariski closure $\G$ of $\Ga$ in
$\GL(n, \CC)$ is reductive. As $b^1(\Ga)=0$ one can assume that $\G$
is simple.

\begin{lemma} \label{projective}
The group $\Ga$ is virtually the fundamental group of a smooth projective
variety.
\end{lemma}

\begin{proof}
This follows from \cite[theor.1]{zuo1} generalizing \cite[Main
Theorem]{mok}, and the faithfulness of $\rho$, as follows.

By \cite[theor1.(b)]{zuo1} there exists an
analytic morphism with connected fibers $l : X\lo M$, $M$ normal
irreducible projective variety, such that
$\rho$ factorizes through $\pi_1(M)$. Let $\sigma_M : \hat{M} \lo M$ be a resolution of
singularities of $M$. One can find a modification $\sigma_X : \hat{X}
\lo X$ and a map $\hat{l} : \hat{X} \lo \hat{M}$ with connected fibers
such that the diagram
$$
\xymatrix{
\hat{X} \ar[r]^{\hat{l}} \ar[d]_{\sigma_{X}} & \hat{M} \ar[d]^{\sigma_{M}} \\
X \ar[r]_l & M}$$ 
commutes.
Taking fundamental groups one obtain the commutative diagram of
groups~:
\begin{equation} \label{proj}
\xymatrix{
\pi_1(\hat{X}) \ar[r]^{\hat{l}_*} \ar[d]_{(\sigma_{X})_*} & \pi_1(\hat{M}) \ar[d]^{(\sigma_{M})_*} \\
\pi_1(X) \ar[r]_{l_{*}} & \pi_1(M)} \qquad .
\end{equation}

As $\sigma$ is a modification and $X$ is smooth the induced map
$(\sigma_X)_{*}: \pi_1(\hat{X}) \lo \pi_1(X) = \Ga $ is an isomorphism. As
$\hat{l}$ and $l$ are connected the morphisms $\hat{l}_*$ and $l_*$
are surjective. As $\rho$ is faithful and factorizes through $\pi_1(M)$ 
the diagram~(\ref{proj}) implies that $\Ga = \pi_1(\hat{M})$ and the result.
\end{proof}

Replacing $X$ by $\hat{M}$ one can assume that $X$ is a smooth connected
projective complex variety of complex dimension $n$.
If $n=1$ the variety $X$ is a smooth curve, of genus $\geq
1$ as $\Ga$ is infinite. In particular $X$ is a $K(\Ga, 1)$ and $b^2(\Ga) = b^2(X) =1$ thus
the conclusion of theorem~\ref{linear} is valid in this case. Thus we
can assume $n\geq 2$. 

Once more one can assume that $X$ is as surface ($n=2$) by replacing
$X$ by a sufficiently ample smooth complete intersection surface in $X$.
By theorem~\ref{eyssi1} we have a diagram of Shafarevich maps 
$$ 
\xymatrix{
\tilde{X}/ \ker \rho \ar[d] \ar[r]^{\widetilde{\sh_\rho}} &
\widetilde{S_\rho}(X) \ar[d]_{\pi} \\
X \ar[r]_{\sh_{\rho}} & \sh_\rho(X)}\;\;.
$$ 
Moreover by theorem~\ref{eyssi2} the normal space
$\widetilde{S_\rho}(X)$ is Stein.

\begin{lemma}
The map $\rho : \Ga \lo \G((\OO_{K})_S)$ factorizes through
$\pi_1(\sh_{\rho}(X))$.
\end{lemma}

\begin{proof}
By lemma~\ref{leEyssi} the harmonic map $h$ factorizes (equivariantly) through $\widetilde{S}_{\rho}(X)$~:
$$ 
\xymatrix{
\ti{X}/\ker \rho \ar[r]^h \ar[d]_{\widetilde{\sh_{\rho}}}  &R_S \\
\widetilde{S}_{\rho}(X) \ar@{-->}[ru] &}\;\;.$$
The group $\Ga$  acts properly discontinuously on $R_S $. By replacing
$\Ga$ if necessary by a finite index subgroup we can assume that $\Ga$
acts freely. Then the previous diagram induces a factorization~:
$$ 
\xymatrix{
X  \ar[r]^h \ar[d]_{\sh_{\rho}}  &R_S/\Ga \\
\sh_{\rho}(X) \ar@{-->}[ru] & }\;\;.$$
At the level of $\pi_1$ it gives a diagram
\begin{equation} \label{e8}
\xymatrix{
\Ga \ar[r]^{h_*=\text{Id}} \ar[d]_{{\sh_{\rho}}_{*}}  &\Ga \\
\pi_1(\sh_{\rho}(X)) \ar@{-->}[ru] & }\;\;,
\end{equation}
and the result.
\end{proof}

\begin{lemma}
Under the assumptions of theorem~\ref{linear} the map $\sh_\rho : X \lo \sh_\rho(X)$ is a modification
and the morphism $(\sh_\rho)_* : \Ga \lo \pi_1(\sh_\rho(X))$ is an isomorphism.
\end{lemma}
\begin{proof}
By definition the map $\sh_\rho$ is surjective thus $\sh_\rho(X)$ is of (complex) dimension at most $2$. By the previous lemma the
representation $\rho$ factorises through
$\pi_1(\widetilde{S}_{\rho})$. As $\Ga= \rho(\Ga)$ is infinite the group
$\pi_1(\widetilde{S}_{\rho})$ too. Thus $\sh_{\rho}(X)$ is not a point
nor $\proj^1 \CC$. If $\sh_\rho(X)$ is  a (smooth) curve of genus $g
\geq 1$ then $b^1(X) \geq b^1(\sh_\rho(X)) >0$, contradiction to
equality~(\ref{e7}). Finally $\sh_\rho(X)$ is of dimension $2$. By
definition the fibers of $\sh_\rho$ are connected thus $\sh_\rho : X
\lo \sh_\rho(X)$ is a modification.

As $\sh_\rho : X \lo \sh_\rho(X)$ is a modification and $\sh_\rho(X)$
is normal we obtain that the map ${\sh_{\rho}}_* : \Ga \lo
\pi_1(\sh_\rho(X))$ is surjective. By considering the
diagram~(\ref{e8}) we obtain that this map is an isomorphism. 
\end{proof}

As $p: X \lo M:=\sh_\rho(X)$ is a modification inducing an isomorphism on fundamental groups the
universal cover $\ti{X}$ is also a modification of $\ti{M}$. Thus the
map $ \ti{p}: \ti{X} \lo \ti{M}$ can be realized topologically as a
cofibration with cofiber $C_{\tilde{p}}$ a union of suspension of complex curves. As
$\ti{M}$ is a Stein surface (possibly singular) it has the homotopy
type of a $2$-dimensional $CW$-complex. 
Writing the (co)homology long exact sequence for the pair $(M,X)$ we
obtain that $X$ does not have any (co)homology in degree larger than
$2$. Thus $X$ still satisfies the hypotheses
of theorem~\ref{top}, which implies the result.

\diam

%%%%%%%%%%%%%%%%%%%%%%%%%%%%%%%%%%
\section{Application to $p$-adic lattices}  \label{padic}

\subsection{Proof of theorem~\ref{p-adic}}

More generally one can prove the following
\begin{theor} \label{p-adic2}
Let $G =  \prod_{v \in S} G_v$, where $S$ denotes a finite set and
$G_v$, $v \in S$, is the group of $F_v$-points of an algebraic
group $\G_v$, with reductive neutral component, over a non-archimedean
local field $F_v$ of characteristic $0$. 
Suppose that $d(G):= \sum_{v \in S}\text{rk}_{F_{v}} \G_v$ satisfies
$d(G)\geq 6$.

Then an irreducible cocompact lattice $\Ga \subset G$ is never K\"ahler.
\end{theor} 
\begin{proof}

Clearly the group $\Ga$ is linear. In particular $\Ga$ admits a finite index
subgroup which is torsion-free. As a group is K\"ahler if and only if
any finite index subgroup of $\Ga$ is K\"ahler, we can assume that $\Ga$ is torsion-free.

By theorem~\cite[theor. 6.2]{BorSer} $\Ga$ is a duality
group of dimension $d(G)$. 

Moreover by a famous result of Garland \cite{gar} (under some
restriction on the residual characteristic of the $F_v$'s, extended by
Casselman \cite{cas} to the general case) any finite index subgroup
$\Ga'$ of $\Ga$ satisfies 
$$H^i(\Ga', \CC) = 0 \quad \text{if} \; i \not = 0, d(G). $$

Thus theorem~\ref{p-adic} immediately follows from
theorem~\ref{linear}.

\begin{rem}
Garland's theorem is proven in the case of a single $v$ but the proof easily
extends to the case under consideration.
\end{rem}

\end{proof}

\subsection{Theorem~\ref{p-adic2} and Simpson's integrality conjecture}
Notice that the conclusion of theorem~\ref{p-adic2} is predicted for
$d(G) \geq 2$ by Simpson's integrality conjecture \cite{si1} saying that any rigid representation of a
K\"ahler group should be integral. 

Indeed for simplicity let us assume that $S = \{v\}$ and $\G_v$ is of adjoint type. Fix an
embedding $F_v \stackrel{\sigma}{\hookrightarrow} \CC$ and consider the natural faithful
representation $\rho: \Ga \lo \G_v^{\sigma} (\CC)$ (where $\G_v^\sigma
:= \G_v \times_{F_{v}, \sigma} \CC$). As $H^1(\Ga, \Ad \rho) =0$ the
representation $\rho$ is rigid. 

If $\Ga$ were a K\"ahler group then by
Simpson's integrality conjecture there would exist a number field $K
\stackrel{\tau_{0}}{\hookrightarrow} \CC$ and a $K$-form
$\G$ of $\G_v^{\sigma} $ such that $\rho(\Ga) \subset
\tau_0(\G(\mathcal{O}_K))$. 

By Margulis's superrigidity theorem \cite{mar}, for any $\tau$ in
the finite set $S_{\infty}$ of Archimedean places of $K$ the
representation $\tau \circ \rho : \Ga \lo \G^\tau(\CC)$ has bounded
image~: there exists a maximal compact subgroup $U_\tau$ of
$\G^\tau(\CC)$ such that $\tau(\rho(\Ga)) \subset U_\tau$. Thus $$\rho
(\Ga) \subset \G(\mathcal{O}_K)  \cap \prod_{\tau \in
  S_{\infty}} U_\tau\;\;.$$

However the intersection $\G(\mathcal{O}_K)  \cap \prod_{\tau \in
  S_{\infty}} U_\tau$ is finite as $\G(\mathcal{O}_K) $ is discrete
in $\prod_{\tau \in  S_{\infty}} \G^\tau(\CC)$.
Contradiction to the fact that the lattice $\Ga$ of $G_v$ is infinite.

%%%%%%%%%%%%%%%%%%%%%%%%%%%%%%%%%%
\appendix
\section{What does conjecture~\ref{conj1} mean~?} \label{explanation}

\subsection{In terms of group extensions}
Recall that for $\Ga$ a group and $A$ an Abelian group with trivial
$\Gamma$-module structure the group $H^2(\Ga, A)$ classifies the
central $A$-extensions
$$ 0 \lo A \lo \tilde{\Gamma} \lo \Gamma \lo 1$$
of $\Ga$. As $b_2(\Ga) = \textnormal{rk} H^2(\Ga, \ZZ)$ the
conjecture~\ref{conj1} means that any infinite K\"ahler group admits
(after maybe passing to a finite index subgroup) a non-trivial
$\ZZ$-central extension (which does not trivialize when restricted to
any finite index subgroup).

\subsection{In topological terms}
For any group $\Ga$ the universal coefficients exact
sequence yields the isomorphism
$$H^2(\Ga, \RR) = \Hom_\RR(H_2(\Ga, \RR), \RR)\;\;.$$
In particular it is equivalent to show $b_2(\Ga) >0$ or $b^2(\Ga)>0$.
For any reasonable topological space $M$ with fundamental group $\Ga$
the universal cover $\tilde{M}$ is a principal $\Ga$-cover of
$M$. Thus it defines (uniquely in the homotopy category) a morphism $c: M \lo B\Ga$
from $M$ to the classifying Eilenberg-MacLane space $B\Ga =
K(\Ga,1)$. The induced morphism 
$c_{*}:H_*(M, \R) \lo H_*(\Ga, \R)$ is easily seen to be an
isomorphism in degree $1$ and an epimorphism in degree $2$~: $$H_2(M, \R) \twoheadrightarrow H_2(\Ga, \R)\;\;.$$
Dualy~:$$H^2(\Ga, \R) \hookrightarrow H^2(M, \R) \;\;.$$
How can we caracterize the quotient $H_2(\Ga, \R)$ of $H_2(M, \R)$~?
In fact this quotient first appeared in Hopf's work on the Hurewicz
morphism comparing homotopy and homology~:
\begin{theor} [Hopf] \label{Hopf}
Let $N$ be a paracompact topological space. Let $c: N \lo B\pi_1(N)$
be the classifying morphism and $h: \pi_*(N) \lo H_*(N, \ZZ)$ the
classical Hurewicz morphism.
Then the sequence of Abelian groups
\begin{equation} \label{equ1}
\pi_2(N) \stackrel{h}{\lo} H_2(N, \ZZ) \stackrel{c_*}{\lo} H_2(\pi_1(N),
\ZZ) \lo 0
\end{equation} is exact.
\end{theor} 

Cohomologically~:
\begin{corol} \label{corol3}
Let $N$ be a paracompact topological space and $\pi_2(N)
\otimes_\ZZ \R \stackrel{h}{\lo} H_2(N, \R)$ the Hurewicz
morphism. Then~:
$$H^2(\pi_1(N), \R) = \{[\om] \in H^2(N, \R)\; /\; \forall \phi:S^2
\lo N, \; <[\om], \phi_* [S^2]>=0\} \subset H^2(N, \R)\;\;,$$ where 
$<\cdot, \cdot>: H^2(N, \R) \times H_2(N, \R) \lo \R$ is the
natural non-degenerate pairing between homology and cohomology.
\end{corol}

\begin{rem}
Nowadays theorem~\ref{Hopf} is a direct
application of the Leray-Cartan spectral sequence.
\end{rem}

Carlson-Toledo's conjecture can thus be restated~:
\begin{conj} \label{conj2}
Let $\Ga$ be an infinite K\"ahler group. There exists a compact
K\"ahler manifold $M$ with $\pi_1(M)$ a finite index subgroup of $\Ga$
such that the Hurewicz morphism $\pi_2(M)\otimes_\ZZ
\RR \lo H_2(M, \RR)$ is not surjective.
\end{conj}

A stronger conjecture is then~:
\begin{conj} \label{conj3}
Let $M$ be a compact K\"ahler manifold with infinite fundamental group
$\Ga$. There exists a finite \'etale cover $M'$ of $M$ such that the Hurewicz morphism $\pi_2(M)\otimes_\ZZ
\RR \lo H_2(M', \RR)$ is not surjective.
\end{conj}

Of course $\pi_2(M)$ is nothing else than $H_2(\tilde{M}, \ZZ)$ where
$\tilde{M}$ denotes the universal cover of $M$. The cohomological
version of the previous conjecture gives~:
\begin{conj} \label{conj4}
Let $M$ be a compact K\"ahler manifold with infinite fundamental group
$\Ga$. Then the natural map 
$$\varinjlim_{M'} H^2(M', \RR) \lo H^2(\tilde{M}, \RR)$$ is not
injective (where the injective limit is taken over the projective
system of \'etale finite cover of $M$). 
\end{conj}

Notice that for any compact manifold $M$ and any finite \'etale cover
$M'$ of $M$ the arrow $H^2(M, \RR) \lo H^2(M', \RR)$ is injective
by the projection formula.

\subsection{In terms of $\CC^*$-bundles}
Recall that for a reasonable topological space $M$ the group $H^2(M,
\ZZ)$ canonically identifies with the group 
$\mathcal{L}(M)$ of principal $\CC^*$-bundles~: on the one hand
$H^2(M, \ZZ) = [M, K(\ZZ, 2)]$, on the other hand 
$\mathcal{L}(M)= [M, B\CC^*]$. But both $K(\ZZ, 2)$ and $B\CC^*$ have
as canonical model the infinite projective space $\CC
\proj^{\infty}$. Thus conjecture~\ref{conj1} states that that any
infinite K\"ahler group $\Ga$ admits a finite index subgroup $\Ga'$
whose classifying space $B\Ga'$ supports a non-trivial
$\CC^*$-torsor. Let $M$ be a compact K\"ahler manifold 
with infinite fundamental group $\Ga$. As $\tilde{M}$ is the homotopy
fiber of $M \lo B\Gamma$ the conjecture~\ref{conj4}
says there exists a finite \'etale cover $M'$ of $M$ and a non-trivial
$\CC^*$-torsor $L'$ on $M'$ whose pull-back to $\tilde{M}$ becomes
trivial. In these statements ``non-trivial'' means ``with non-trivial
{\em rational} first Chern class''.

\subsection{Equivalence} \label{equivalence}
The equivalence of these $3$ points of view is clear. Given $M$ and a
$\CC^*$-torsor $L$ on $M$ one can consider the long homotopy exact
sequence for the fibration $L \lo M$~:
\begin{equation} \label{long}
 \cdots \lo \pi_2(L) \lo \pi_2(M) \stackrel{c_1}{\lo}\pi_1(\CC^*) = \ZZ \lo
\pi_1(L) \lo \pi_1(M) \lo 1 \;\;.
\end{equation}
The boundary map $c_1: \pi_2(M) \lo \ZZ$ is just the first Chern class
map of $L$ restricted to $\pi_2(M)$~: if $[\alpha] \in \pi_2(M)$ is
represented by $\alpha: S^2 \lo M$ then $c_1([\alpha])=
<\alpha^* (c_1(L)), [S^2]> \in \ZZ$. As $H_2(M, \ZZ) = \pi_2(M)$ and
$H^2(\tilde{M}, \RR)$ is dual to $H_2(M, \RR)$ the torsor
$p^*(L)$ is trivial if and only if $c_1:
\pi_2(M) \lo \RR$ is zero. Then the long exact sequence~(\ref{long})
gives the short exact sequence~: 
$$ 1 \lo \ZZ \lo \pi_1(L) \lo \pi_1(M) \lo 1 \;\;.$$

%%%%%%%%%%%%%%%%%%
\section{Carlson-Toledo's conjecture and rigidity} \label{sec2}
The strongest evidence for Carlson-Toledo's conjecture is the
following folkloric theorem, essentially due to Lefschetz, Simpson and
Reznikov (a proof is provided in \cite{KlKoMau})~:

\begin{theor} \label{prop2}
Let $\Ga$ be a K\"ahler group. If $\Ga$ admits a linear representation $\rho:
\Ga \lo G$, with $G$ the linear group $GL(V)$ of a finite dimensional
complex vector space $V$ or the isometry group of a Hilbert space, satisfying
$H^1(\Ga, \rho) \not = 0$,
then $b_2(\Ga) >0$.
\end{theor}

Recall the following definitions~:
\begin{defi}  \label{rigi}
A finitely generated group $\Ga$ is said
\begin{itemize}
\item[(1)] {\em rigid} if for any  $n \in \NN$
$\Ga$ admits only finitely many {\em simple} linear repr\'esentations
of dimension $n$ (up to conjugacy).
\item[(2)] {\em reductive} if all its linear finite-dimensional
  complex representations are semi-simple.
\item[(3)]{\em schematically rigid} if for any $n \in \NN$
and any representation $\rho: \Ga \lo \GL(n, \CC)$ one has $H^1(\Ga,
\Ad \rho)=0$.
\end{itemize}
\end{defi}

\begin{rem}
In \cite{lubma} a rigid group is called $S$-rigid 
($SS$ for semi-simple).
\end{rem}

Let $A(\Ga)$ be the pro-(affine)-algebraic completion of $\Ga$. For $n
\in \NN$ let $A_n(\Ga)$ be the quotient
$A(\Ga)/K_n(\Ga)$ where $K_n(\Ga)$ denotes the intersection of the
kernels of all representations $A(\Ga \lo GL(n, \CC)$. All
$\Ga$-representations of dimension $n$ factorizes uniquely through
$A_n(\Ga)$. One easily shows that $\Ga$ is
rigid if and only if for all $n \in \NN$ the group $A_n(\Ga)$
is an affine algebraic group (i.e. of finite dimension)
(c.f. \cite[theorem A]{blmm}).

\sspace
In the definition~\ref{rigi} each condition implies the previous one~:

As the $\GL(n, \CC)$-orbit of a semi-simple representation $\rho \in \Hom(\Ga,
\GL(n, \CC)$ is closed condition $(2)$ implies condition~$(1)$. On the
other hand there exists non-reductive rigid groups~: for example
$\SL(n, \ZZ) \ltimes \ZZ^n$, $n \geq 3$.

Any schematically-rigid group is rigid of course. As the $GL(n,
\CC)$-orbit in $\Hom(\Ga, \GL(n, \CC)$ of any representation contains
a semi-simple representation in its closure any schematically-rigid
group is reductive. To say that $\Ga$ is schematically rigid is
equivalent to saying that 
$A(\Ga)$ is pro-reductive. There exist rigid groups which are not
schematically rigid~: c.f. \cite[section
5]{blmm}.

As a corollary of theorem~\ref{prop2} one obtains~:

\begin{corol}
Let $\Ga$ an infinite K\"ahler group. If $\Ga$ does not satisfy
Carlson-Toledo's conjecture then necessarily~:
\begin{itemize}
\item[(a)] $\Ga$ has Kazhdan's property $(T)$.
\item[(b)] $\Ga$ is schematically rigid. 
\end{itemize}
\end{corol}

%%%%%%%%%%%%%%%%%%%%%%%%%%%%%%%%%%%%%%%%%%%%%%%%%%%%%%%%%%%%%%%%%%%%%%%%%%%

\lspace
\noi
{Bruno Klingler

\noi Institut de Math{\'e}matiques de Jussieu, Paris, France

\noi
e-mail~: klingler@math.jussieu.fr}

\end{document}